\documentclass[journal,  comsoc, onecolumn]{IEEEtran}

\usepackage[T1]{fontenc}

\usepackage{amssymb}
\usepackage{amsmath}

\usepackage{enumitem} 
\usepackage{tikz}
\usepackage{tkz-graph}
\usepackage{graphicx}
\usepackage{balance}

\ifCLASSINFOpdf
\else
\fi
\usepackage{graphicx}

\usepackage{amsmath}
\usepackage{color}

\usepackage[hidelinks]{hyperref}
\usepackage{xcolor}
\hypersetup{
	colorlinks,
	linkcolor={red!80!black},
	citecolor={blue!80!black},
	urlcolor={blue!80!black}
}

\usepackage{url}

\usepackage{mathrsfs}
\newtheorem{thm}{Theorem}
\newtheorem{lem}[thm]{Lemma}

\newtheorem{prob}{Problem}

\newtheorem{defn}[thm]{Definition}

\newtheorem{rem}[thm]{Remark}

\begin{document}
%
\title{$\mathcal{H}_2$ and $\mathcal{H}_\infty$ Suboptimal Distributed Filter Design \\ for Linear Systems}
%
%
%
\author{Junjie~Jiao,~Harry~L.~Trentelman,~\IEEEmembership{Fellow,~IEEE,} and M.~Kanat~Camlibel,~\IEEEmembership{Member,~IEEE}


\thanks{The authors are with the Bernoulli Institute for Mathematics, Computer Science and Artificial Intelligence, University of Groningen, Groningen, 9700 AV, The Netherlands (Email: \href{mailto:j.jiao@rug.nl}{j.jiao@rug.nl};
\href{mailto:h.l.trentelman@rug.nl}{h.l.trentelman@rug.nl}; \href{mailto:m.k.camlibel@rug.nl}{m.k.camlibel@rug.nl}).}
}
%
%
%

\markboth{Journal of \LaTeX\ Class Files,~Vol.~14, No.~8, \today} 
{Shell \MakeLowercase{\textit{et al.}}: Bare Demo of IEEEtran.cls for IEEE Communications Society Journals}
%



\maketitle

\begin{abstract}
This paper investigates the $\mathcal{H}_2$ and $\mathcal{H}_\infty$  suboptimal  distributed filtering problems for continuous time linear systems.
Consider a  linear system monitored by a number of  filters, where each of the filters receives only part of the measured output of the system.
Each filter can communicate with the other filters according to an a priori given strongly connected weighted directed graph.
The aim is to design filter gains that guarantee the $\mathcal{H}_2$ or  $\mathcal{H}_\infty$ norm of the transfer matrix from  the disturbance input to the output estimation error to be smaller than an a priori given upper bound, while all local filters reconstruct the full system state asymptotically.
We provide a centralized design method for obtaining such  $\mathcal{H}_2$ and $\mathcal{H}_\infty$ suboptimal distributed filters.
The proposed design method is illustrated by a simulation example.
\end{abstract}


%
\IEEEpeerreviewmaketitle

\section{Introduction}\label{sec_intro}
\IEEEPARstart{R}{ecent} years have witnessed an increasing interest in problems of state estimation for spatially constrained large-scale systems. 
Such problems are relevant  in  applications, such as power grids \cite{PowerGrid2012spm}, industrial plants \cite{plantProcess2003} and wireless sensor networks \cite{Olfati-Saber2005-2}. 
Due to the physical constraints, the measured output of these systems is often monitored by a sensor network, consisting of a number of local sensors. 
Each of these local sensors makes use of its local measurements and then communicates with the other local sensors. 
In this way, all of these  sensors together are able to estimate the state of the system asymptotically.
In this problem setting,   one of the main challenges is that none of the local sensors by itself is able to estimate the system state by using its own local measurements.
Consequently,  standard estimation methods do not directly apply anymore.

The distributed estimation problem has  been  mainly studied in two research directions, namely, distributed observer design and distributed Kalman filtering.
%
In \cite{Martins2017}, an augmented state observer was proposed to cast the distributed observer design problem into a decentralized control problem for linear systems, using the notion of `fixed modes' \cite{Davison1973}.
Later on, in \cite{Morse2018}, the results in \cite{Martins2017} were extended and a more general form of distributed observers was provided, allowing the rate of convergence  of the observer to be freely assignable.
In  \cite{morse2017cdc}, for time-varying communication graphs,  a hybrid observer was introduced to distributedly estimate the state of a linear system.
Based on observability decompositions, the problem of distributed observer design was also investigated in \cite{han_observer}, \cite{Sundaram2018}  and \cite{shim2016cdc}.
In \cite{MITRA2019108487}, an attack resilient algorithm was introduced to address the distributed estimation problem when certain nodes are compromised by adversaries.

On the other hand, much attention in the literature has also been devoted to distributed filtering problems.
A Kalman-filter-based distributed filter was proposed in \cite{Olfati-Saber2005-1},  \cite{Olfati-Saber2009} and  \cite{Olfati-Saber2007}. There, the proposed methods employ a two-step strategy: a state update rule based on a Kalman-filter and a data fusion step based on consensus.
%
In \cite{UGRINOVSKII20111}, a distributed robust filtering problem was addressed using dissipativity theory. 
%
Later on in \cite{UGRINOVSKII2013160}, the results of \cite{UGRINOVSKII20111} were generalized  to the case that the communication graph is allowed to randomly change.
Recently, in \cite{Shim2016}, a distributed Kalman-Bucy filtering problem was studied, using the idea of averaging the dynamics of heterogeneous multi-agent systems \cite{Shim2016heterogeneous}.

%
Different from the existing work, in the present paper,  we will consider two deterministic versions of the distributed optimal filtering problem  for linear systems, i.e., the  $\mathcal{H}_2$ and $\mathcal{H}_\infty$ distributed filtering problems.
%
%
Given a linear system and a network of local filters, each local filter receives a portion of the  measured output  of the system and then exchanges its  state with that of its neighboring  local filters. 
Together, these local filters form a distributed filter. 
We introduce $\mathcal{H}_2$ and $\mathcal{H}_\infty$ performances to quantify the influence of the disturbances  on the output estimation error. The distributed optimal filtering problem  is then to find suitable filter gain matrices such that the associated $\mathcal{H}_2$ or $\mathcal{H}_\infty$ performance  is minimized, while the states of all local filters asymptotically track the system state.
However, due to {\em non-convexity}, this problem is difficult to solve in general.
Therefore, in this paper we will address a {\em suboptimality} version of this problem. 
The objective of the present paper is then  to design suitable filter gain matrices such that the $\mathcal{H}_2$ or $\mathcal{H}_\infty$  performance  is smaller than an a priori given tolerance. 
%
The main contributions of this paper are the following:
\begin{enumerate}
	\item We establish conditions for the existence of suitable filter gains in terms of solvability of LMI's for both   the $\mathcal{H}_2$ and $\mathcal{H}_\infty$ suboptimal distributed filtering problem. For the $\mathcal{H}_2$
	filtering problem, all except one of these LMI's will always turn out to be solvable.
	\item We provide conceptual algorithms for obtaining suitable  $\mathcal{H}_2$ and $\mathcal{H}_\infty$ suboptimal distributed filters, respectively.
\end{enumerate}

This paper is organized as follows. 
In Section \ref{sec_pre}, we review some basic results on graph theory, detectability properties of linear systems,   and the $\mathcal{H}_2$ and $\mathcal{H}_\infty$ performance of linear systems. 
Subsequently, in Section \ref{sec_prob} we formulate the   $\mathcal{H}_2$ and $\mathcal{H}_\infty$ suboptimal distributed filtering problems. 
We then provide design methods for obtaining such distributed filters in Section \ref{sec_filterdesign}. 
In Section \ref{sec_simulation} we provide a simulation example to illustrate our design method.
Finally, in Section \ref{sec_conclusion} we formulate our conclusions.

\subsection*{Notation}
We denote by  $\mathbb{R}$ the field of real numbers and by  $\mathbb{R}^n$ the space of $n$ dimensional vectors over $\mathbb{R}$.
We write $\mathbf{1}_N$ for the $n$ dimensional column vector with all its  entries equal to $1$.
%
For a given matrix $A$, we write $A^\top$ to denote  its transpose and $A^{-1}$ its inverse (if exists).
For a symmetric matrix $P$, we denote $P>0$ if it  is positive definite and $P < 0$  if its negative definite.
We denote the identity matrix of dimension $n \times n$ by $I_n$.
A matrix is called Hurwitz if all its eigenvalues have negative real parts. 
The trace of a square matrix $A$ is denoted by $\textnormal{tr}(A)$.
We denote by $\text{diag}(d_1, d_2, \ldots, d_n)$ the $n \times n$ diagonal matrix with $d_1, d_2,\ldots, d_n$ on the diagonal.
Given matrices $R_i \in \mathbb{R}^{m \times m}$, $i=1,2,\ldots, n$, we denote by $\text{blockdiag}(R_i)$ the $nm \times nm$ block diagonal matrix with $R_1, R_2,\ldots, R_n$ on the diagonal and we denote by $\text{col}(R_i)$ the $nm \times m$ column block matrix $\left(R_1^\top, R_2^\top,\ldots, R_n^\top\right)^\top$.
The Kronecker product of two matrices $A$ and $B$ is denoted by $A \otimes B$.
For a linear map $A: \mathcal{X} \to \mathcal{Y}$, the kernel and image of $A$ are denoted by $\ker (A) : = \{ x \in \mathcal{X} \mid A x =0 \}$ and ${\rm im}(A) := \{ Ax \mid x \in \mathcal{X} \}$, respectively.

\section{Preliminaries}\label{sec_pre}

\subsection{Graph Theory}\label{subsec_graph}
A weighted directed graph is denoted by $\mathcal{G} = (\mathcal{V}, \mathcal{E}, \mathcal{A})$, where $\mathcal{V} = \{ 1,2,\ldots, N \}$ is the finite nonempty node set, $\mathcal{E} \subset \mathcal{V} \times \mathcal{V}$ is the edge set of ordered pairs $(i,j)$  and $\mathcal{A} = [a_{ij}]$ is the associated adjacency matrix with nonnegative entries.
The entry $a_{ji}$ of the adjacency matrix $\mathcal{A}$ is the weight associated with the edge $(i,j)$ and $a_{ji}$ is nonzero if and only if $(i,j) \in \mathcal{E}$. 
%
Given a graph $\mathcal{G}$, a directed path from node $1$ to node $p$ is a sequence of edges $(k, {k+1})$, $k = 1,2,\ldots, p-1$.
A graph is called {\em strongly connected} if for any pair of distinct nodes $i$ and $j$, there exists a directed path from $i$ to $j$.
A graph is called simple if $a_{ii} =0$, i.e., the graph does not contain self-loops.
A graph is called  undirected if $(i,j) \in \mathcal{E}$ implies   $(j,i) \in \mathcal{E}$.
A simple undirected graph is called {\em connected} if for each pair of nodes $i$ and $j$ there exists a path from $i$ to $j$.

Given a graph $\mathcal{G}$, the degree matrix of $\mathcal{G}$ is denoted by $\mathcal{D} = \textnormal{diag}(d_1,d_2,\ldots, d_N)$ with $d_i = \sum_{j=1}^N a_{ij}$.
The Laplacian matrix of  $\mathcal{G}$ is defined as $L := \mathcal{D} - \mathcal{A}$.
If $\mathcal{G}$ is a weighted directed graph, the associated Laplacian matrix $L$ has a zero eigenvalue corresponding to the eigenvector $\mathbf{1}_N$. 
If moreover $\mathcal{G}$  is strongly connected,  then all the other eigenvalues lie in the open right half-plane.
%

For strongly connected weighted directed graphs, we review the following lemma \cite{MeiJie2016}:
\begin{lem}\label{lem_laplacian}
	Let $\mathcal{G}$ be a strongly connected weighted directed graph with Laplacian matrix $L$.
	Then there exists a unique row vector 
	$\theta = \left(\theta_1, \theta_2,\ldots, \theta_N\right) $,
	where $\theta_1, \theta_2, \ldots, \theta_N$ are all positive real numbers, 	such that 
	$\theta  L = 0$ and $\theta  \mathbf{1}_N = N$.
	%
	Define 
	$	\Theta := \textnormal{diag} (\theta_1, \theta_2,\ldots, \theta_N),$
	then the matrix 
	$	\mathcal{L} := \Theta L +L^\top \Theta$
	is a positive semi-definite matrix associated with a connected weighted undirected  graph. 
	%
\end{lem}

\subsection{Detectability and Detectability Decomposition}\label{subsec_detect}
In this subsection, we  review  detectability and the detectability decomposition of linear systems.
Consider the  linear  system
\begin{equation}\label{sys}
\begin{aligned}
\dot{x} &= {A} x,\\
y &= Cx ,
\end{aligned}
\end{equation}
where $x \in \mathbb{R}^n$ represents the state and $y \in \mathbb{R}^p$ the measured output. 
The matrices $A$ and $C$ are of suitable dimensions.

Let $p(s)$ be the characteristic polynomial of $A$.
Then $p(s)$ can be factorized as 
\[
p(s) = p_-(s) p_+(s),
\]
where $p_-(s)$ and $p_+(s)$ have roots in the open left half-plane and the closed right half-plane, respectively.
The undetectable subspace of the pair $(C, A)$ is defined as
\begin{equation*}
\mathcal{S} := \mathcal{N} \cap \ker\big( p_+(A)\big),
\end{equation*}
where
\begin{equation*}
\mathcal{N} := \ker 
\begin{pmatrix}
C \\
CA \\
\vdots \\
CA^{n-1}
\end{pmatrix}.
\end{equation*}
The pair $(C, A)$ is detectable if and only if  
$
\mathcal{S} = \{ 0 \}
$, see e.g. \cite{Wonham1985book}.

There exists an orthogonal matrix  $T \in \mathbb{R}^{n \times n}$  such that the pair $(C, A)$ is transformed into the detectability decomposition form
\begin{equation*}
T^\top A T = 
\begin{pmatrix}
A_{11} & 0 \\
A_{21} & A_{22}
\end{pmatrix}, 
\quad
C T = 
\begin{pmatrix}
C_{1} & 0
\end{pmatrix},
\end{equation*}
where $A_{11} \in \mathbb{R}^{v \times v}$, $A_{21} \in \mathbb{R}^{(n-v) \times v}$, $A_{22} \in \mathbb{R}^{(n-v)\times (n-v)}$,  $C_{1} \in \mathbb{R}^{p \times v}$ and
the pair $(C_1, A_{11})$ is detectable.
In addition, if we partition $T = (T_1\ T_2)$, where $T_1$ contains the first $v$ columns, then the undetectable subspace is given by 
\[
\textnormal{im}(T_2) = \mathcal{S}.
\]
Since $T$ is orthogonal, we also have
\[
\textnormal{im}(T_1) = \mathcal{S}^\perp.
\]


\subsection{$\mathcal{H}_2$ and $\mathcal{H}_\infty$ Performance of Linear Systems}\label{subsec_linear_analysis}
In this subsection, we review the $\mathcal{H}_2$ and $\mathcal{H}_\infty$ performance of a linear system with external disturbances.
Consider the  linear  system
\begin{equation}\label{sys_d}
\begin{aligned}
\dot{x} &= {A} x + {E} d,\\
y &= {C}x ,
\end{aligned}
\end{equation}
where $x \in \mathbb{R}^n$ is the state, $d \in \mathbb{R}^q$ the external disturbance and $y \in \mathbb{R}^p$ the measured output. 
The matrices ${A}$, ${C}$ and ${E}$ are of suitable dimensions.

We first review the  $\mathcal{H}_2$ performance of the system \eqref{sys_d}. 
Let
$
T_d(t) = {C} e^{{A}t}{E}
$
be the impulse response of \eqref{sys_d}.
Then the associated $\mathcal{H}_2$ performance is  defined to be the square of its $L_2$-norm, given by
\begin{equation}\label{cost_d}
J = \int_{0}^{\infty} \textnormal{tr} \left[ T_d^{\top}(t) T_d(t) \right] dt.
\end{equation}
Note that the performance \eqref{cost_d}  is finite if the system \eqref{sys_d} is internally stable, i.e., ${A}$ is Hurwitz.

The following well-known result provides a necessary and sufficient condition under which \eqref{sys_d} is internally stable and  (\ref{cost_d}) is smaller than a given upper bound (see  e.g. \cite{IWASAKI1994421}, \cite{Scherer1997}):

\begin{lem}\label{sys_d_lem}
	%
	Let  $\gamma > 0$. 
	Then the system \eqref{sys_d} is internally stable  and $J < \gamma$ if and only if there exists $P >0$ satisfying
	\begin{align*}
	{A}^{\top} P + P {A}  + {C}^{\top} {C} &< 0,  \\
	\textnormal{tr}\left( {E}^{\top} P {E} \right) &<\gamma.  
	\end{align*}
\end{lem}

Next, we review the  $\mathcal{H}_\infty$ performance of   the system \eqref{sys_d}. 
Let $T_d(s) = {C} (sI_n -{A})^{-1}{E}$ be the transfer matrix of \eqref{sys_d}.
If $A$ is Hurwitz, then the  $\mathcal{H}_\infty$ performance of \eqref{sys_d} is  defined as the $\mathcal{H}_\infty$ norm of $T_d(s)$, given by
\begin{equation}\label{norm_infty}
||T_d ||_\infty := \sup_{\omega \in \mathbb{R}} \sigma(T(j\omega)),
\end{equation}
where $\sigma(T_d(j\omega))$ is the maximum singular value of the complex matrix $T_d(j\omega)$.

The  well-known bounded real lemma  provides a necessary and sufficient condition under which  \eqref{sys_d} is stable and (\ref{norm_infty}) is smaller than a given upper bound (see  e.g. \cite{kemin_zhou_book}, \cite{harry_book}):

\begin{lem}\label{hinf_lem}
	%
	Let  $\gamma > 0$. 
	Then  the system \eqref{sys_d} is internally stable  and $||T_d ||_\infty < \gamma$ if and only if there exists $P >0$ such that
	\begin{equation*}
	{A}^{\top} P + P {A}  + \frac{1}{\gamma^2}P {E} {E}^{\top}P + {C}^\top {C} < 0.
	\end{equation*}
\end{lem}

In the next section, we will formulate the $\mathcal{H}_2$ and  $\mathcal{H}_\infty$ distributed filter design problems that  will be addressed in this paper.

\section{Problem Formulation}\label{sec_prob}
Consider the  finite-dimensional linear time-invariant system
\begin{equation}\label{sys_general}
\begin{aligned}
\dot{x} & = A x +   E d,\\
y &= C x +D d, \\
z &= H x ,
\end{aligned}
\end{equation}
where $x \in \mathbb{R}^n$ is the  state, $d\in \mathbb{R}^q$   the external disturbance, $y \in \mathbb{R}^r$  the measured output and $z\in \mathbb{R}^p$ the output  to be estimated.
The matrices $A$, $C$, $D$, $E$ and $H$  are  of suitable dimensions.
%

The standard  optimal filtering problem for the system \eqref{sys_general} is to find a filter that takes $y$ as input and returns an optimal estimate $\zeta$ of $z$, while the filter state asymptotically tracks the state $x$ of \eqref{sys_general}. 
Here, `optimal' means that the $\mathcal{H}_2$ or $\mathcal{H}_\infty$ norm of the transfer matrix from $d$ to the estimation error $z -\zeta$ is minimized over all such filters.
In that problem setting, however, a standing assumption is that one single filter is able to acquire the complete    measured output $y$ of the system.

In the present paper,  we relax this assumption. More specifically, we assume that the measured output $y$ of  \eqref{sys_general} is not available to one single filter, but is observed by $N$  local filters. 
Moreover, each local filter only acquires a certain portion of the measured output, namely,
\[
y_i = C_i x + D_i d, 
\]
where $y_i \in \mathbb{R}^{r_i}$, $C_i \in \mathbb{R}^{r_i \times n}$ and $D_i \in \mathbb{R}^{r_i \times q}$, for $i = 1,2,\ldots , N$.
Here, the matrices $C_i$ and $D_i$ are obtained by partitioning 
\begin{equation*}
C = \begin{pmatrix}
C_1 \\
C_2 \\
\vdots\\
C_N
\end{pmatrix},\quad
D = \begin{pmatrix}
D_1 \\
D_2 \\
\vdots\\
D_N
\end{pmatrix}.
\end{equation*}
Clearly, the original output $y$ of \eqref{sys_general} has then been partitioned as
\[
y = \begin{pmatrix}
y_1 \\
y_2 \\
\vdots \\
y_N	
\end{pmatrix}
\]
and  $\sum_{i=1}^{N} r_i = r$. 
%
In this paper it will be a standing assumption that the pair $(C, A)$ is detectable.
We will also assume that none of the pairs $(C_i, A)$ is detectable itself.
If, for at least one $i$, the pair $(C_i, A)$ is detectable, 
the distributed filtering problem boils down to the standard optimal filtering problem.

In our distributed case, each local filter makes use of   the portion of the measured output that it acquires and will then communicate with its neighboring local filters by exchanging filter state information. 
In this way, the local filters  will together form a distributed filter.
Following \cite{han_observer} and \cite{shim2016cdc}, we propose a distributed filter of the form
\begin{equation}\label{dis_observer}
\begin{aligned}
\dot{w}_i &= A w_i + G_i (y_i - C_i w_i) + F_i  \sum_{j =1}^N a_{ij} (w_j -w_i) ,\\
\zeta_i & = Hw_i, \quad i = 1,2,\ldots , N,
\end{aligned}
\end{equation}
where $w_i \in \mathbb{R}^n$ is the state of the $i$th local filter and $\zeta_i \in \mathbb{R}^p$ is the associated output. The matrices $G_i \in \mathbb{R}^{n \times r_i}$ and $F_i \in \mathbb{R}^{n \times n}$ are local filter gains to be designed. The coefficients $a_{ij}$ are the entries of the adjacency matrix $\mathcal{A}$ of the communication graph. In this paper, it will be a standing assumption that this graph is a strongly connected weighted directed graph.

%
%

%

For the $i$th local filter, we introduce the associated local state estimation error $e_i$ and local output estimation error $\eta_i$ as 
\begin{align*}
e_i &:= x- w_i , \\
\eta_{i}&:= z - \zeta_i, \quad i = 1,2,\ldots, N.  
\end{align*}
The dynamics of the $i$th local error system is then given by
\begin{equation*}\label{local_error}
\begin{aligned}
\dot{e_{i}} &= (A -G_i C_i) e_{i} + F_i \sum_{j =1}^N a_{ij} (e_{j} -e_{i}) + (E- G_i D_i) d ,\\
\eta_{i} &= H e_{i}, \quad i = 1,2,\ldots , N.
\end{aligned}
\end{equation*}
Denote 
$e  = (e_1^\top, e_2^\top, \ldots, e_N^\top)^\top$, $\eta = (\eta_1^\top, \eta_2^\top, \ldots, \eta_N^\top)^\top$ 
and
\begin{align*}
\bar{A} &:= \text{blockdiag}( A -G_i C_i) \in \mathbb{R}^{nN \times nN}, \\
\bar{F} &:= \text{blockdiag}(F_i) \in \mathbb{R}^{nN \times nN}, \\
\bar{E} &:= \text{col}(E -  G_i D_i) \in \mathbb{R}^{nN \times q}.
\end{align*}
The global error system is then given by
\begin{equation}\label{global_error_sys}
\begin{aligned}
\dot{e} &= \left(\bar{A} -  \bar{F} (L \otimes I_n)\right) e + \bar{E} d , \\
\eta &= (I_N \otimes H) e,
\end{aligned}
\end{equation}
where $L \in \mathbb{R}^{N \times N}$ is the Laplacian matrix of the communication graph.
The impulse response of the system \eqref{global_error_sys} from the disturbance $d$ to the  output estimation error $\eta$ is equal to
\begin{equation*}
T_d(t) =  (I_N \otimes H) e^{(\bar{A} - \bar{F} (L \otimes I_n))t} \bar{E}.
\end{equation*}
We introduce the global $\mathcal{H}_2$ cost functional 
\begin{equation}\label{cost_h2}
J =  \int_{0}^{\infty} \text{tr} \left[ T_d^\top (t) T_d(t) \right]  dt.
\end{equation}
The  $\mathcal{H}_2$ optimal distributed filtering problem is then the problem of minimizing the  $\mathcal{H}_2$ cost functional \eqref{cost_h2} over all distributed filters \eqref{dis_observer} such that the global error system \eqref{global_error_sys} is internally stable. Note that \eqref{cost_h2} is a function of the local gain matrices $F_1, F_2, \ldots, F_N$ and $G_1, G_2, \ldots, G_N$.

Unfortunately, due to the particular form of  \eqref{dis_observer}, this optimization problem is, in general, non-convex and it is unclear whether a closed-form solution exists. 
Therefore, instead of trying to find an {\em optimal }solution,  we will address a version of this problem that only requires {\em suboptimality}.
More concretely, we aim at designing a distributed filter such that the error system \eqref{global_error_sys} is internally stable and the  $\mathcal{H}_2$ performance \eqref{cost_h2} is smaller than an a priori given tolerance $\gamma$. 
In that case, we say that the distributed filter \eqref{dis_observer} is   $\mathcal{H}_2$ $\gamma$-suboptimal:
\begin{defn}\label{def1}
	Let $\gamma >0$. The distributed filter \eqref{dis_observer} is called $\mathcal{H}_2$ $\gamma$-suboptimal if:
	\begin{enumerate}
		\item  for all $i = 1,2, \ldots, N$, whenever $d=0$,  we have that
		$\lim_{t \to \infty}  \big(x(t) - w_i(t)\big) \to 0$ for all initial conditions on \eqref{sys_general} and \eqref{dis_observer}.
		\item  the associated performance \eqref{cost_h2} satisfies $J < \gamma$. 
	\end{enumerate} 
\end{defn}

Correspondingly, the  $\mathcal{H}_2$ suboptimal distributed filtering problem that we will address  is the following:
\begin{prob}\label{Prob1}
	%
	%
	Let $\gamma >0$.
	For $i = 1,2, \ldots, N$, find gain matrices  $G_i \in \mathbb{R}^{n \times r_i}$ and $F_i \in \mathbb{R}^{n \times n}$   such that the distributed filter \eqref{dis_observer} is $\mathcal{H}_2$ $\gamma$-suboptimal.
\end{prob}

In addition to the distributed filtering problem with $\mathcal{H}_2$ performance, in this paper we will also consider the version of this problem with $\mathcal{H}_\infty$ performance. Obviously, the transfer matrix of the system  \eqref{global_error_sys} from the disturbance $d$ to the output estimation error $\eta$ is equal to
\[
T_d(s) = (I_N \otimes H) \left(sI_{nN} - (\bar{A} -  \bar{F} (L \otimes I_n)) \right)^{-1}\bar{E}.
\]
The  $\mathcal{H}_\infty$ performance of the distributed filter \eqref{dis_observer} is given by the $\mathcal{H}_\infty$  norm $||T_d||_\infty$ of $T_d(s)$.
The problem that we will then consider is to design a distributed filter \eqref{dis_observer} such that the error system \eqref{global_error_sys} is internally stable and its  $\mathcal{H}_\infty$ performance is smaller than an a priori given tolerance $\gamma$. 
In that case, we say that the distributed filter \eqref{dis_observer} is   $\mathcal{H}_\infty$ $\gamma$-suboptimal:
\begin{defn}\label{def2}
	Let $\gamma >0$.
	The distributed filter \eqref{dis_observer}  is called $\mathcal{H}_\infty$ $\gamma$-suboptimal if:
	\begin{enumerate}
		\item  for all $i = 1,2, \ldots, N$, whenever $d=0$,  we have that
		$\lim_{t \to \infty}  \big(x(t) - w_i(t)\big) \to 0$ for all initial conditions on \eqref{sys_general} and \eqref{dis_observer}.
		\item $||T_d ||_\infty < \gamma$.
	\end{enumerate}
\end{defn}

Correspondingly,  the  $\mathcal{H}_\infty$ suboptimal distributed filtering problem that we will address  is the following:
\begin{prob}\label{Prob2}
	%
	%
	Let $\gamma >0$.
	For $i = 1,2, \ldots, N$, find gain matrices  $G_i \in \mathbb{R}^{n \times r_i}$ and $F_i \in \mathbb{R}^{n \times n}$  such that the distributed filter \eqref{dis_observer} is  $\mathcal{H}_\infty$ $\gamma$-suboptimal.
\end{prob}

\section{$\mathcal{H}_2$ and $\mathcal{H}_\infty$ Suboptimal  Distributed Filter Design}\label{sec_filterdesign}
In this section, we will address Problems \ref{Prob1} and \ref{Prob2} introduced above and provide design methods for obtaining suboptimal distributed filters. 

As we have explained before,  the $i$th local filter \eqref{dis_observer} receives only a certain portion of the measured output, namely,
\[
y_i= C_i x + D_i d,\quad i = 1,2, \ldots, N.
\]
%
In order to proceed, we first apply   orthogonal transformations to the pairs $(C_i, A)$. 
For $i = 1,2, \ldots, N$, let $T_i$ be an orthogonal matrix such that the pair $(C_{i}, A)$ is transformed into the  detectability decomposition form
\begin{equation}\label{Detectability_AC}
T_i^\top A T_i = 
\begin{pmatrix}
A_{i11} & 0 \\
A_{i21} & A_{i22}
\end{pmatrix},\quad
C_i T_i = 
\begin{pmatrix}
C_{i1} & 0
\end{pmatrix} ,
\end{equation}
where   $A_{i11} \in \mathbb{R}^{v _i\times v_i}$, $A_{i21} \in \mathbb{R}^{(n-v_i) \times v_i}$, $A_{i22} \in \mathbb{R}^{(n-v_i)\times (n-v_i)}$,  $C_{i1} \in \mathbb{R}^{r_i \times v_i}$ and the pair $(C_{i1}, A_{i11})$ is detectable. The integer $v_i$ is equal to the dimension of the othogonal complement of the undetectable subspace of the pair $(C_i,A)$.
Accordingly, partition
\begin{equation}\label{Detectability_EH}
T_i^\top E  = 
\begin{pmatrix}
E_{i1} \\
E_{i2}
\end{pmatrix},\quad
H T_i = 
\begin{pmatrix}
H_{i1} & H_{i2}
\end{pmatrix},
\end{equation}
where $E_{i1} \in \mathbb{R}^{v _i\times q}$, $E_{i2} \in \mathbb{R}^{(n-v _i)\times q}$, $H_{i1} \in \mathbb{R}^{p \times v _i}$ and $H_{i1} \in \mathbb{R}^{p \times (n - v _i)}$.

Using the fact that $(C_{i1}, A_{i11})$ is detectable, let $Q_{i1}$ be any positive definite solution to 
\begin{equation}\label{AR_ineq_Q_directed}
A_{i11} Q_{i1} + Q_{i1} A_{i11}^{\top} - Q_{i1} C_{i1}^{\top} C_{i1} Q_{i1}  < 0.
\end{equation}
Then, by defining 
\begin{equation}\label{G_i1}
G_{i1}  := Q_{i1} C_{i1}^\top,
\end{equation}
the matrix $A_{i11} - G_{i1}C_{i1}$ is Hurwitz.

In the sequel, we will make use of the  transformed matrices \eqref{Detectability_AC} and \eqref{Detectability_EH}  and the gain matrix \eqref{G_i1} to obtain filter gains that solve Problems \ref{Prob1} and  \ref{Prob2}.
Before presenting the main results of this paper, we will first provide a lemma that will be essential for later use. This lemma is a generalization of \cite[Lemma 4]{han_observer}, and connects the Laplacian matrix of the communication graph with detectability properties of the system \eqref{sys_general}. 

\begin{lem}\label{lem_L_directed}
	Let $\mathcal{L} := \Theta L +L^\top \Theta$, where $\Theta$ is  defined as in Lemma \ref{lem_laplacian}.
	Define $T := \textnormal{blockdiag}( T_i) \in \mathbb{R}^{nN \times nN}$, where the $T_i$ are the orthogonal matrices introduced in \eqref{Detectability_AC} and \eqref{Detectability_EH}.
	Let $m_i>0$ and
	\begin{equation*}
	M_i := 
	\begin{pmatrix}
	m_i I_{v_i }& 0 \\
	0 & 0_{n-v_i}
	\end{pmatrix}, \quad i = 1,2, \ldots, N.
	\end{equation*} 
	Define $M := \textnormal{blockdiag}(M_i)$. Then, 
	\begin{equation}\label{ineq_L_directed}
	T^\top (\mathcal{L} \otimes I_n) T +  M > 0.
	\end{equation}
\end{lem}
The proof of Lemma \ref{lem_L_directed} can be given by adapting the proof of \cite[Lemma 4]{han_observer}, replacing the observability decomposition by the detectability decomposition. We omit the details here.

In the next two subsections, we will deal with the design of $\mathcal{H}_2$  and $\mathcal{H}_\infty$ suboptimal   distributed filters, respectively.

\subsection{$\mathcal{H}_2$ Suboptimal Distributed Filter Design}\label{subsec_H2}

In this subsection, we will  provide a design method for obtaining  $\mathcal{H}_2$ suboptimal distributed  filters.
More specifically, we aim at finding a distributed filter such that the global error system \eqref{global_error_sys} is stable and the associated $\mathcal{H}_2$ performance \eqref{cost_h2} is less than an a priori given tolerance.

The next lemma expresses  the existence of suitable gain matrices $F_i$ and $G_i$,  $i = 1,2,\ldots, N$
in terms of solvability of LMI's.

\begin{lem} \label{lem_Prob1_directed}
	%
	%
	Let $\gamma >0$. 
	Let  the matrices  $T$, $M$ and $\mathcal{L}$ be as introduced in Lemma \ref{lem_L_directed}.
	Let $\epsilon>0$ be such that 
	\begin{equation}\label{epsilon_h2}
	T^\top (\mathcal{L} \otimes I_n) T +  M > \epsilon I_{nN}.
	\end{equation}
	Let  $G_{i1}$ be as defined in \eqref{G_i1}.
	For $i = 1,2,\ldots,N$, assume there exist $\kappa >0$, $P_{i1} >0$ and $P_{i2} >0$ satisfying
	\begin{equation}  \label{ineq_perform1_directed}
	\begin{pmatrix}
	\Phi_i + H_{i1}^\top H_{i1} + \kappa(m_i -\epsilon) I_{v_i}& A_{i21}^\top P_{i2} + H_{i1}^\top H_{i2}  \\
	P_{i2} A_{i21} + H_{i2}^\top H_{i1} & \Psi_i
	\end{pmatrix} <0 
	\end{equation}
	and 
	\begin{equation}
	\sum_{i=1}^{N}\textnormal{tr} \left[ (E_{i1} - G_{i1}D_i)^\top  P_{i1} (E_{i1} - G_{i1}D_i)  + E_{i2}^\top P_{i2} E_{i2} \right] < \gamma, \label{ineq_gamma1_directed} 
	\end{equation}
	where  
	\begin{align}
	\Phi_i & := A_{i11}^\top P_{i1} + P_{i1} A_{i11} - C_{i1}^\top  G_{i1}^\top P_{i1} - P_{i1}  G_{i1} C_{i1} , \label{Phi}\\
	\Psi_i & := P_{i2} A_{i22} + A_{i22}^\top P_{i2} +H_{i2}^\top H_{i2} - \kappa \epsilon I_{n - v_i}.\label{Psi}
	\end{align}
	For $i = 1,2, \ldots, N$, define gain matrices $F_i$ and $G_i$ by
	\begin{equation}\label{F_i}
	F_i := \kappa \theta_i 
	T_i
	\begin{pmatrix}
	P_{i1}^{-1} & 0 \\
	0 & P_{i2}^{-1} 
	\end{pmatrix}
	T_i^\top 
	\end{equation}
	and
	\begin{equation}\label{G_i}
	G_i := T_i
	\begin{pmatrix}
	G_{i1} \\
	0
	\end{pmatrix}.
	\end{equation}
	Then the corresponding distributed filter \eqref{dis_observer}  is $\mathcal{H}_2$ $\gamma$-suboptimal.
\end{lem}

\begin{IEEEproof}
	First, it follows from  \eqref{ineq_L_directed} in Lemma \ref{lem_L_directed} that there exists   $\epsilon>0$ such that \eqref{epsilon_h2} holds.
	Next, note that \eqref{ineq_perform1_directed} is equivalent to 
	\begin{equation}\label{ineq1}
	\begin{aligned}
	\textnormal{blockdiag}\begin{pmatrix}
	\Phi_i + H_{i1}^\top H_{i1}& A_{i21}^\top P_{i2}  + H_{i1}^\top H_{i2} \\
	P_{i2} A_{i21} + H_{i2}^\top H_{i1} & P_{i2} A_{i22} + A_{i22}^\top P_{i2} +H_{i2}^\top H_{i2}
	\end{pmatrix} &\\
	+ \kappa (M -\epsilon I_{nN})
	&<0.
	\end{aligned}
	\end{equation}
	Using \eqref{epsilon_h2}, it follows from \eqref{ineq1} that 
	\begin{equation}\label{ineq2}
	\begin{aligned}
	\textnormal{blockdiag}\begin{pmatrix}
	\Phi_i + H_{i1}^\top H_{i1}& A_{i21}^\top P_{i2}  + H_{i1}^\top H_{i2} \\
	P_{i2} A_{i21} + H_{i2}^\top H_{i1} & P_{i2} A_{i22} + A_{i22}^\top P_{i2} +H_{i2}^\top H_{i2}
	\end{pmatrix} &\\
	- \kappa T^\top (\mathcal{L} \otimes I_n) T
	&<0.
	\end{aligned}
	\end{equation}
	Let
	\begin{equation}\label{P}
	P := \textnormal{blockdiag}(P_i),\quad
	P_i := 
	T_i
	\begin{pmatrix}
	P_{i1} & 0 \\
	0 & P_{i2}
	\end{pmatrix}
	T_i^\top.
	\end{equation}
	%
	Clearly, $P>0$.
	By using  \eqref{F_i}, \eqref{G_i}, \eqref{P}, \eqref{Detectability_AC} and \eqref{Detectability_EH}, then \eqref{ineq2} holds if and only if
	\begin{equation}\label{ineq3}
	\begin{aligned}
	\bar{A}^\top P  + P \bar{A}  -  (L^\top \otimes I_n) \bar{F}^\top P + P \bar{F}(L \otimes I_n) 
	+ I_N \otimes H^\top H <0
	\end{aligned}
	\end{equation}
	holds,
	where $\bar{F} := \textnormal{blockdiag}(F_i)$ and $F_i$ is defined by \eqref{F_i}.
	Therefore, there exist $\kappa >0$, $P_{i1} >0$ and $P_{i2} >0$ such that  \eqref{ineq_perform1_directed} holds for $i=1,2,\ldots,N$ if and only if there exists $P>0$ of the form \eqref{P} such that \eqref{ineq3} holds.
	Since the solutions of \eqref{ineq_perform1_directed} also satisfy \eqref{ineq_gamma1_directed}, we obtain 
	\begin{equation}\label{trace}
	\textnormal{tr} \left( \bar{E}^\top P \bar{E} \right) < \gamma.
	\end{equation}
	Finally, since \eqref{ineq3} and \eqref{trace} have a solution $P>0$,  it follows from Lemma \ref{sys_d_lem} that the error system \eqref{global_error_sys} is internally stable and  $J < \gamma$.
	Thus the distributed filter \eqref{dis_observer} with \eqref{G_i} and \eqref{F_i} is $\mathcal{H}_2$ $\gamma$-suboptimal.
\end{IEEEproof}

\begin{rem}\label{Remark1}
	In Lemma \ref{lem_Prob1_directed},  the choice of the parameters $m_i>0$ is arbitrary. The parameter $\epsilon>0$ should be chosen sufficiently small so that \eqref{epsilon_h2} holds. The gain $G_{i}$ is defined by \eqref{G_i}.
	Then, of course,  the question arises: for chosen $m_i>0$, $\epsilon >0$ and $G_{i}$, how can we find the smallest $\gamma >0$ such that the corresponding distributed filter \eqref{dis_observer} is $\mathcal{H}_2$ $\gamma$-suboptimal? This requires to find the smallest $\gamma$ such that the LMI's \eqref{ineq_perform1_directed} and \eqref{ineq_gamma1_directed} are solvable. It is well known that this can be done by using a standard bisection algorithm, see e.g. \cite[page 115]{kemin_zhou_book}.
\end{rem}

\begin{rem} 
	Lemma \ref{lem_Prob1_directed} states that if there exist solutions $\kappa >0$, $P_{i1} >0$ and $P_{i2} >0$ satisfying \eqref{ineq_perform1_directed} and \eqref{ineq_gamma1_directed}, then the  distributed filter \eqref{dis_observer} with gain matrices \eqref{F_i} and \eqref{G_i} is $\mathcal{H}_2$ $\gamma$-suboptimal. 
	There, the inequality  \eqref{ineq_gamma1_directed} is a global condition  for checking suboptimality.
	In fact, such suboptimality condition can also be  checked locally.
	Indeed, if  for $i =1,2,\ldots,N$ there exist  solutions satisfying \eqref{ineq_perform1_directed} and 
	\begin{equation*}
	\textnormal{tr} \left[ (E_{i1} - G_{i1}D_i)^\top  P_{i1} (E_{i1} - G_{i1}D_i)  + E_{i2}^\top P_{i2} E_{i2} \right] < \frac{\gamma}{N},
	\end{equation*}
	then the corresponding  distributed filter \eqref{dis_observer} with \eqref{F_i} and \eqref{G_i} is $\mathcal{H}_2$ $\gamma$-suboptimal. 
\end{rem}
Lemma \ref{lem_Prob1_directed} provides a condition for the existence of suitable gain matrices $F_i$ and $G_i$ in terms of solvability of LMI's. 
In the next theorem, we show that, in fact, the LMI's \eqref{ineq_perform1_directed} in Lemma \ref{lem_Prob1_directed} always have solutions. In fact, we can take $P_{i2}$ to be the identity matrix of dimension  $n - v_i$ and $P_{i1}$ to be the unique solution of a given Lyapunov equation. In this way we obtain the following conceptual algorithm for computing suitable gain matrices.
\begin{thm}\label{main_1}
	%
	Let $\gamma >0$.
	Then an  $\mathcal{H}_2$ $\gamma$-suboptimal  distributed filter of the form \eqref{dis_observer}  is obtained as follows:
	\begin{enumerate}
		\item\label{step_theta} Compute $\theta = \left(\theta_1, \theta_2,\ldots, \theta_N\right)$ with $\theta_i 		>0$  such that $\theta  L = 0$ and  $\theta  \mathbf{1}_N = N$.
	\end{enumerate}
	Then, for $i =1,2,\ldots, N$:
	\begin{enumerate}
		\setcounter{enumi}{1}
		\item  Compute orthogonal matrices $T_i$ that put $A,$ $ C_i,$ $E$  and $H$ into the form \eqref{Detectability_AC} and \eqref{Detectability_EH}.
		\item\label{epsilon} Take $m_i =1$ and compute $\epsilon >0$ such that
		\begin{equation}\label{epsi}
		T^\top (\mathcal{L} \otimes I_n) T +  M > \epsilon I_{nN}.
		\end{equation}
		\item Compute  $Q_{i1}>0$ satisfying \eqref{AR_ineq_Q_directed}. Define $G_{i1} : = Q_{i1} C_{i1}^\top$.
		\item \label{lmi_kappa}  Take $\kappa >0$ sufficiently large such that 
		\begin{equation}\label{kappa}
		\begin{aligned}
		& A_{i22} + A_{i22}^\top  +H_{i2}^\top H_{i2} - \kappa \epsilon I_{n - v_i} \\
		&\quad + \frac{1}{\kappa \epsilon}(A_{i21} + H_{i2}^\top H_{i1})(A_{i21} + H_{i2}^\top H_{i1})^\top <0.
		\end{aligned}
		\end{equation}
		\item Compute $P_{i1} > 0$ satisfying the Lyapunov equation
		\begin{equation}\label{lyap}
		\begin{aligned}
		(A_{i11} - G_{i1} C_{i1})^\top P_{i1} + P_{i1} (A_{i11} - G_{i1} C_{i1}) \\
		\qquad + H_{i1}^\top H_{i1}  + \kappa I_{v_i} &= 0.
		\end{aligned}
		\end{equation}	
		\item 	Define gain matrices $F_i$ and $G_i$ by 
		\begin{equation}\label{Gi_Fi}
		F_i  := \kappa \theta_i
		T_i
		\begin{pmatrix}
		P_{i1}^{-1} & 0 \\
		0 & I_{n-v_i}
		\end{pmatrix}
		T_i^\top,\quad
		G_i  := T_i
		\begin{pmatrix}
		G_{i1} \\
		0
		\end{pmatrix}.
		\end{equation} 
	\end{enumerate}
	Then for all $\gamma >0$ satisfying 
	\begin{equation}\label{tolerance}
	\sum_{i=1}^{N}\textnormal{tr} \left[ (E_{i1} - G_{i1}D_i)^\top P_{i1} (E_{i1} - G_{i1}D_i)  + E_{i2}^\top  E_{i2} \right] < \gamma,
	\end{equation}
	the corresponding distributed filter \eqref{dis_observer} with gain matrices \eqref{Gi_Fi} is  $\mathcal{H}_2$ $\gamma$-suboptimal.
	
\end{thm}
\begin{IEEEproof}
	Using Lemma \ref{lem_L_directed}, by choosing $m_i =1$ for $i =1,2,\ldots, N$,  there exists $\epsilon>0$ such that \eqref{epsi} holds.
	Next, for $i =1,2,\ldots, N$, we choose $\kappa >0$ sufficiently large such that \eqref{kappa} holds .
	%
	%
	Since $Q_{i1}$ is a positive definite solution of \eqref{AR_ineq_Q_directed} and $G_{i1} : = Q_{i1} C_{i1}^\top$, then the matrix $A_{i11} - G_{i1} C_{i1}$ is Hurwitz. 
	Consequently, for $i =1,2,\ldots, N$, the Lyapunov equation \eqref{lyap} has  unique solution $P_{i1} >0$.
	Since \eqref{kappa} holds and $-\kappa \epsilon I_{v_i} < 0$,  by using the Schur complement,  we obtain
	\begin{equation}\label{eq6}
	\begin{pmatrix}
	-\kappa \epsilon  I_{v_i}& A_{i21}^\top   + H_{i1}^\top H_{i2}  \\
	A_{i21} + H_{i2}^\top H_{i1} & \tilde{\Psi}_i
	\end{pmatrix} <0, \quad i =1,2,\ldots, N,
	\end{equation}
	where
	$\tilde{\Psi}_i:=  A_{i22} + A_{i22}^\top   +H_{i2}^\top H_{i2} - \kappa \epsilon I_{n - v_i}.$
	Using \eqref{lyap} and $P_{i2} = I_{n - v_i}$, 	it then follows that \eqref{ineq_perform1_directed} holds.
	%
	
	On the other hand, by taking $P_{i2} = I_{n - v_i}$ in \eqref{ineq_gamma1_directed}, we obtain  \eqref{tolerance}.
	It then follows from Lemma \ref{lem_Prob1_directed} that the corresponding distributed filter is $\mathcal{H}_2$ $\gamma$-suboptimal.
\end{IEEEproof}

%


\begin{rem}\label{Remark3}
	%
	Note that, in step \ref{step_theta}) of Theorem \ref{main_1}, we need to compute the left eigenvector $\theta$ of the Laplacian matrix corresponding to the eigenvalue $0$. 
	This requires so-called global information on the communication graph. 
	This dependency on global information can be removed   using algorithms that compute   left eigenvectors of the Laplacian matrix in a distributed fashion, see e.g. \cite{Hadjicostis2016} or \cite{Qu2017}.
	On the other hand, in step \ref{epsilon}) we need to compute $\epsilon$. 
	To do so, we need  knowledge of  the orthogonal matrices $T_i$, the matrix $M$ and the Laplacian matrix $\mathcal{L}$, which is global information.
	Also in step \ref{lmi_kappa}), we  need  to find one $\kappa$ that satisfy  \eqref{kappa} for $i =1,2,\ldots, N$.  
	Note that, however, we can always take $\epsilon>0$ sufficiently small and $\kappa>0$ sufficiently large such that \eqref{epsi} and  \eqref{kappa} hold, respectively. 
	This might however lead to an achievable tolerance $\gamma$ that is very large, giving poor suboptimality of the corresponding distributed filter. 
\end{rem}

In general, the computation of our suboptimal filters requires global information, so  cannot be performed in a decentralized fashion. This is in contrast with the decentralized computation of distributed state observers as described in \cite{8943159}.

\subsection{$\mathcal{H}_\infty$ Suboptimal  Distributed Filter Design}\label{subsec_Hinf}
%
In this subsection, we will  provide a  method for obtaining $\mathcal{H}_{\infty}$ suboptimal distributed filters.
More concretely, we aim at finding, for a given tolerance $\gamma>0$,  a distributed filter such that the global error system \eqref{global_error_sys} is stable and  $||T_d ||_\infty < \gamma$.

The next lemma expresses the existence of suitable gain matrices $F_i$ and $G_i$,  $i = 1,2,\ldots, N$ in terms of solvability of $N$ nonlinear matrix inequalities.

\begin{lem} \label{lem_Prob2}
	%
	%
	Let $\gamma >0$.
	Let  the matrices $T$, $M$  and $\mathcal{L}$ be as introduced in Lemma \ref{lem_L_directed}.
	Let $\epsilon>0$ be such that 
	\begin{equation}\label{epsilon_hinf}
	T^\top (\mathcal{L} \otimes I_n) T +  M > \epsilon I_{nN}.
	\end{equation}
	Let  $G_{i1}$ be as defined  in \eqref{G_i1}.
	For $i = 1,2,\ldots,N$, assume there exist $\kappa >0$, $P_{i1} >0$ and $P_{i2} >0$ satisfying 
	\begin{equation}\label{ineq_perform_inf}
	\begin{pmatrix}
	\Phi_i + \kappa(m_i -\epsilon) I_{v_i}  & \Omega_i \\
	\Omega_i^\top & \Psi_i- \kappa \epsilon I_{n - v_i}
	\end{pmatrix} <0, 
	\end{equation}
	where 
	\begin{align}
	\Phi_i &= (A_{i11} - G_{i1}^\top C_{i1} )^\top P_{i1} + P_{i1}(A_{i11} - G_{i1}^\top C_{i1} ) \nonumber  \\
	&\quad	+ \frac{1}{\gamma^2} P_{i1}(E_{i1}  - G_{i1} D_i)(E_{i1}  - G_{i1} D_i)^\top P_{i1}+ H_{i1}^\top H_{i1} , \label{Phi_inf}\\
	\Omega_i &= A_{i21}^\top P_{i2} + H_{i1}^\top H_{i2} +\frac{1}{\gamma^2} P_{i1} (E_{i1}  - G_{i1} D_i)E_{i2}^\top P_{i2}, \label{Omega_inf} \\
	\Psi_i &= P_{i2} A_{i22} + A_{i22}^\top P_{i2} +\frac{1}{\gamma^2} P_{i2}E_{i2} E_{i2}^\top P_{i2}  +H_{i2}^\top H_{i2}  . \nonumber
	\end{align}
	For $i =1,2,\ldots,N$, define gain matrices $F_i$  and $G_i$ by
	\begin{equation}\label{F_i_inf}
	F_i := \kappa \theta_i 
	T_i
	\begin{pmatrix}
	P_{i1}^{-1} & 0 \\
	0 & P_{i2}^{-1} 
	\end{pmatrix}
	T_i^\top
	\end{equation}
	and 
	\begin{equation}\label{G_i_inf}
	G_i := T_i
	\begin{pmatrix}
	G_{i1} \\
	0
	\end{pmatrix} .
	\end{equation}
	Then the corresponding distributed filter \eqref{dis_observer}  is $\mathcal{H}_\infty$ $\gamma$-suboptimal.
\end{lem}

\begin{IEEEproof}
	First, it follows from  \eqref{ineq_L_directed} in Lemma \ref{lem_L_directed} that there exists   $\epsilon>0$ such that \eqref{epsilon_h2} holds.
	Next, note that \eqref{ineq_perform_inf} is equivalent to 
	\begin{equation}\label{ineq_inf1}
	\begin{aligned}
	\textnormal{blockdiag} 	\begin{pmatrix}
	\Phi_i  & \Omega_i \\
	\Omega_i^\top & \Psi_i  
	\end{pmatrix} 
	+ \kappa (M -\epsilon I_{nN}) <0  .
	\end{aligned}
	\end{equation}
	Using  \eqref{epsilon_h2}, it then follows from \eqref{ineq_inf1} that 
	\begin{equation}\label{ineq_inf2}
	\begin{aligned}
	\textnormal{blockdiag} 	\begin{pmatrix}
	\Phi_i  & \Omega_i \\
	\Omega_i^\top & \Psi_i  
	\end{pmatrix} 
	- \kappa T^\top (\mathcal{L} \otimes I_n) T  <0  .
	\end{aligned}
	\end{equation}
	Let
	\begin{equation}\label{P_inf}
	P := \textnormal{blockdiag}(P_i),\quad
	P_i := 
	T_i
	\begin{pmatrix}
	P_{i1} & 0 \\
	0 & P_{i2}
	\end{pmatrix}
	T_i^\top.
	\end{equation}
	Clearly, $P>0$.
	By using \eqref{F_i_inf}, \eqref{G_i_inf},  \eqref{P_inf}, \eqref{Detectability_AC} and \eqref{Detectability_EH}, then \eqref{ineq_inf2} holds  if and only if
	\begin{equation}\label{ineq_inf3}
	\begin{aligned}
	\bar{A}^\top P  + P \bar{A} - (L^\top \otimes I_n) \bar{F}^\top P - P \bar{F}(L \otimes I_n)& \\
	+\frac{1}{\gamma^2} P \bar{E}\bar{E}^\top P   +   I_N \otimes H^\top H  &< 0
	\end{aligned}
	\end{equation}
	holds,
	where $\bar{F} := \textnormal{blockdiag}(F_i)$ and $F_i$ is defined by \eqref{F_i_inf}.
	Therefore, there exist  $\kappa >0$, $P_{i1} >0$ and $P_{i2} >0$  such that  \eqref{ineq_perform_inf} holds for $i=1,2,\ldots,N$ if and only if there exists $P>0$ of the form \eqref{P} such that \eqref{ineq_inf3} holds.
	Finally, since \eqref{ineq_inf3} has a solution $P>0$, it follows from Lemma \ref{hinf_lem} that the error system \eqref{global_error_sys} is internally stable and  $||T_d ||_\infty < \gamma$. 
	Thus the distributed filter \eqref{dis_observer} with  \eqref{G_i_inf} and \eqref{F_i_inf} is $\mathcal{H}_\infty$ $\gamma$-suboptimal.
\end{IEEEproof}

\begin{table*}[t]
	\begin{equation}	\label{lmi}
	\begin{pmatrix}
	(A_{i11} - G_{i1}^\top C_{i1} )^\top P_{i1} + P_{i1}(A_{i11} - G_{i1}^\top C_{i1} ) + \kappa(m_i -\epsilon) I_{v_i}  + H_{i1}^\top H_{i1}  &  A_{i21}^\top P_{i2} + H_{i1}^\top H_{i2} &  P_{i1} (E_{i1}  - G_{i1} D_i)  \\
	P_{i2} A_{i21}  + H_{i2}^\top H_{i1}   &    P_{i2} A_{i22} + A_{i22}^\top P_{i2} & P_{i2} E_{i2}  \\
	(E_{i1}  - G_{i1} D_i)^\top P_{i1} & E_{i2}^\top P_{i2}  & -\gamma^2 I_q
	\end{pmatrix} <0.
	\end{equation}
	\medskip
	\hrule
\end{table*}

Lemma \ref{lem_Prob2} provides a condition for  the existence of suitable gain matrices $F_i$ and $G_i$ in terms of solvability of the nonlinear matrix inequalities \eqref{ineq_perform_inf}.
However, these inequalities  are not LMI's.
However, by using  suitable Schur complements, we can transform the inequalities \eqref{ineq_perform_inf}  into LMI's.
In this way we obtain the following conceptual algorithm for computing suitable gain matrices.
\begin{thm}\label{main_2}
	%
	Let $\gamma >0$.
	Then an $\mathcal{H}_\infty$ suboptimal distributed  filter  of the form \eqref{dis_observer} is obtained as follows:
	\begin{enumerate}
		\item \label{step_theta_inf} Compute $\theta = \left(\theta_1, \theta_2,\ldots, \theta_N\right)$  with $\theta_i >0$ such that $\theta  L = 0$ and $\theta  \mathbf{1}_N = N$.
	\end{enumerate}
	For  $i =1,2,\ldots, N$: 
	\begin{enumerate} \setcounter{enumi}{1}
		\item Compute an orthogonal matrix $T_i$ that  puts $A,$ $ C_i,$ $E$ and $H$ into the form \eqref{Detectability_AC} and \eqref{Detectability_EH}.
		\item \label{step_epsilon_inf}  Take arbitrary $m_i >0$  and compute $\epsilon >0$ such that
		\begin{equation}\label{epsi_inf}
		T^\top (\mathcal{L} \otimes I_n) T +  M > \epsilon I_{nN}.
		\end{equation}
		\item Compute  $Q_{i1}>0$ satisfying \eqref{AR_ineq_Q_directed}. Define $G_{i1}  := Q_{i1} C_{i1}^\top$.
		\item  Compute $P_{i1} >0$, $P_{i2} >0$  and $\kappa >0$ such that the inequality \eqref{lmi} (see next page) holds.
		\item Define gain matrices $F_i$  and  $G_i$ by
		\begin{equation}\label{F_iG_i_inf_Thm}
		F_i  := \kappa \theta_i
		T_i
		\begin{pmatrix}
		P_{i1}^{-1} & 0 \\
		0 & I_{n-v_i}
		\end{pmatrix}
		T_i^\top, \quad
		G_i  := T_i
		\begin{pmatrix}
		G_{i1} \\
		0
		\end{pmatrix} .
		\end{equation}
	\end{enumerate}
	Then the corresponding distributed filter \eqref{dis_observer} is $\mathcal{H}_\infty$ $\gamma$-suboptimal.
\end{thm}
\begin{IEEEproof}
	By  taking the appropriate Schur complements in \eqref{lmi}, it follows that \eqref{lmi} hold if and only if  \eqref{ineq_perform_inf} hold. 
	The rest  follows from Lemma \ref{lem_Prob2}.
\end{IEEEproof}

We conclude this section by noting that remarks similar to Remark \ref{Remark1} and Remark \ref{Remark3} hold in the $\mathcal{H}_\infty$ case.

\section{Simulation Example}\label{sec_simulation}
In this section, we will use a simulation example borrowed from \cite{shim2016cdc} to illustrate the conceptual algorithm  in Theorem \ref{main_1} for designing  $\mathcal{H}_2$ suboptimal distributed filters.
Consider the linear time-invariant system
\begin{equation}\label{sys_sim}
\begin{aligned}
\dot{x} & = A x +   E d,\\
y &= C x +D d, \\
z &= H x ,
\end{aligned}
\end{equation}
where 
	\begin{align*}
	A &= 
	\begin{pmatrix}
	0   &  1 &    0 &    0\\
	-1   &  0   &  0   &  0\\
	0    & 0   &  0 &    2\\
	0    & 0 &   -2  &   0
	\end{pmatrix},\quad
	E =
	\begin{pmatrix}
	0.1 \\
	0.1 \\
	0\\
	0.1 
	\end{pmatrix},
	\quad C = I_4,
	\\
	D &= \begin{pmatrix}
	0.1 \\
	0 \\
	0.1\\
	0.1
	\end{pmatrix},\quad
	H = \begin{pmatrix}
	1   &  0    & 0 &     0\\
	0   &  2  &   0   &  0\\
	0   &  0   &  0  &   1
	\end{pmatrix}.
	\end{align*}
The system \eqref{sys_sim} is monitored by four local filters, and each local filter acquires a portion of the measured output $y$, namely,
\begin{equation*}
y_i =C_i 
x + D_id , \quad i = 1,2,3,4,
\end{equation*}
where the matrices $C_i$ and $D_i$ are obtained by partitioning
	\begin{align*}
	C =
	\left(
	\begin{array}{c}
	C_1 \\ \hline
	C_2 \\ \hline
	C_3\\ \hline
	C_4
	\end{array}\right) 
	=
	\left(
	\begin{array}{cccc}
	1 &0 &0& 0 \\ 
	\hline
	0 &1& 0 &0 \\ \hline
	0 &0 &1 &0 \\ \hline
	0 &0& 0 &1
	\end{array}\right) 
	, \quad
	D = 
	\left(
	\begin{array}{c}
	D_1 \\ \hline
	D_2 \\ \hline
	D_3\\ \hline
	D_4
	\end{array}\right) 
	=
	\left(
	\begin{array}{c}
	0.1 \\  \hline
	0 \\ \hline
	0.1\\ \hline
	0.1
	\end{array}\right).
	\end{align*}
The pair $(C,A)$ is detectable but none of the pairs $(C_i, A)$ is detectable.

%
%
%
We assume the four local filters to be of the form \eqref{dis_observer}. 
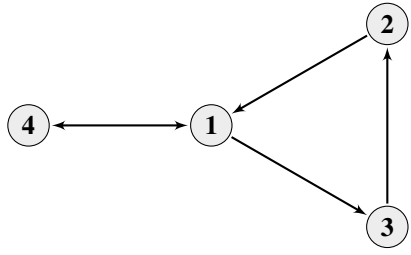
\begin{figure}[t]
	\centering
	\begin{tikzpicture}[scale=0.9]
	\tikzset{VertexStyle1/.style = {shape = circle,
			color=black,
			fill=white!93!black,
			minimum size=0.5cm,
			text = black,
			inner sep = 2pt,
			outer sep = 1pt,
			minimum size = 0.55cm},
		VertexStyle2/.style = {shape = circle,
			color=black,
			fill=black!53!white,
			minimum size=0.5cm,
			text = white,
			inner sep = 2pt,
			outer sep = 1pt,
			minimum size = 0.55cm}
	}
	\node[VertexStyle1,draw](4) at (-3.7,0) {$\bf 4$};
	\node[VertexStyle1,draw](1) at (-1,0) {$\bf 1$};
	\node[VertexStyle1,draw](2) at (1.6,1.5) {$\bf 2$};
	\node[VertexStyle1,draw](3) at (1.6,-1.5) {$\bf 3$};
	\Edge[ style = {<->,> = latex',pos = 0.2},color=black, labelstyle={inner sep=0pt}](4)(1);
	\Edge[ style = {<-,> = latex',pos = 0.2},color=black, labelstyle={inner sep=0pt}](1)(2);
	\Edge[ style = {<-,> = latex',pos = 0.2},color=black, labelstyle={inner sep=0pt}](2)(3);
	\Edge[ style = {<-,> = latex',pos = 0.2},color=black, labelstyle={inner sep=0pt}](3)(1);
	\end{tikzpicture}
	\caption{The communication graph  between the local filters.}
	\label{graph}
\end{figure}
The communication graph between the four local filters is depicted in Figure \ref{graph}. The graph is strongly connected and the associated Laplacian matrix is given by
	\begin{equation*}
	L = \begin{pmatrix}
	2   & -1   &  0  &  -1\\
	0   &  1   & -1   &  0\\
	-1   &  0  &   1  &   0\\
	-1  &   0  &   0   &  1
	\end{pmatrix}.
	\end{equation*}
The normalized left eigenvector $\theta$ of $L$ associated with eigenvalue $0$ is computed to be 
$\theta = \begin{pmatrix}
1 & 1 & 1 & 1
\end{pmatrix}.
$

Next, for $i = 1,2,3,4$, we compute an orthogonal matrix $T_i$ such that the matrices $A,$ $ C_i,$ $E$  and $H$ are transformed into the form \eqref{Detectability_AC} and \eqref{Detectability_EH}.
For  $i =1,2,3,4$, we take $m_i =1$. 
We also compute $\epsilon =0.42$ such that  \eqref{epsi} holds.
Subsequently, for  $i =1,2,3,4$, we solve  \eqref{kappa}  and compute $\kappa = 9.6$.
Following the steps in Theorem \ref{main_1},  gain matrices  $F_i$ and $G_i$ are then computed as
	\begin{equation*}
	\begin{aligned}
	F_1 &=
	\begin{pmatrix}
	0.3636  &  0.0837    &     0     &    0\\
	0.0837    &0.3442    &     0     &    0\\
	0     &    0    &9.6000     &    0\\
	0       &  0    &     0  &  9.6000
	
	\end{pmatrix},\\
	F_2 &= 
	\begin{pmatrix}
	0.3460  & -0.0660    &     0      &   0\\
	-0.0660&    0.3579      &   0   &      0\\
	0      &   0 &   9.6000     &    0\\
	0   &      0   &      0    &9.6000
	\end{pmatrix},
	\end{aligned}
	\end{equation*}
	\begin{equation*}
	\begin{aligned}
	F_3 &= 
	\begin{pmatrix}
	
	9.6000      &   0    &     0     &    0\\
	0 &   9.6000   &      0    &     0\\
	0     &    0  &  0.4274  &  0.0491\\
	0     &    0  &  0.0491  &  0.4217
	\end{pmatrix},\\
	F_4 &= 
	\begin{pmatrix}
	9.6000    &     0     &    0    &     0\\
	0 &   9.6000    &     0    &     0\\
	0   &      0  &  0.4220  & -0.0445\\
	0  &       0  & -0.0445 &   0.4266
	\end{pmatrix},
	\end{aligned}
	\end{equation*}
and
	\begin{equation*}\label{G_sim}
	\begin{aligned}
	G_1 &= 
	\begin{pmatrix}
	0.4445\\
	0.0488\\
	0\\
	0
	\end{pmatrix},\quad
	G_2 = 
	\begin{pmatrix}
	-0.0488\\
	0.4445\\
	0\\
	0
	\end{pmatrix},\\
	G_3 &= 
	\begin{pmatrix}
	0\\
	0\\
	0.4465\\
	0.0248\\
	\end{pmatrix},\quad
	G_4 = 
	\begin{pmatrix}
	0\\
	0\\
	-0.0248\\
	0.4465
	\end{pmatrix}.
	\end{aligned}
	\end{equation*}

\begin{figure}[t]
	\centering
	\includegraphics[width=0.8\columnwidth]{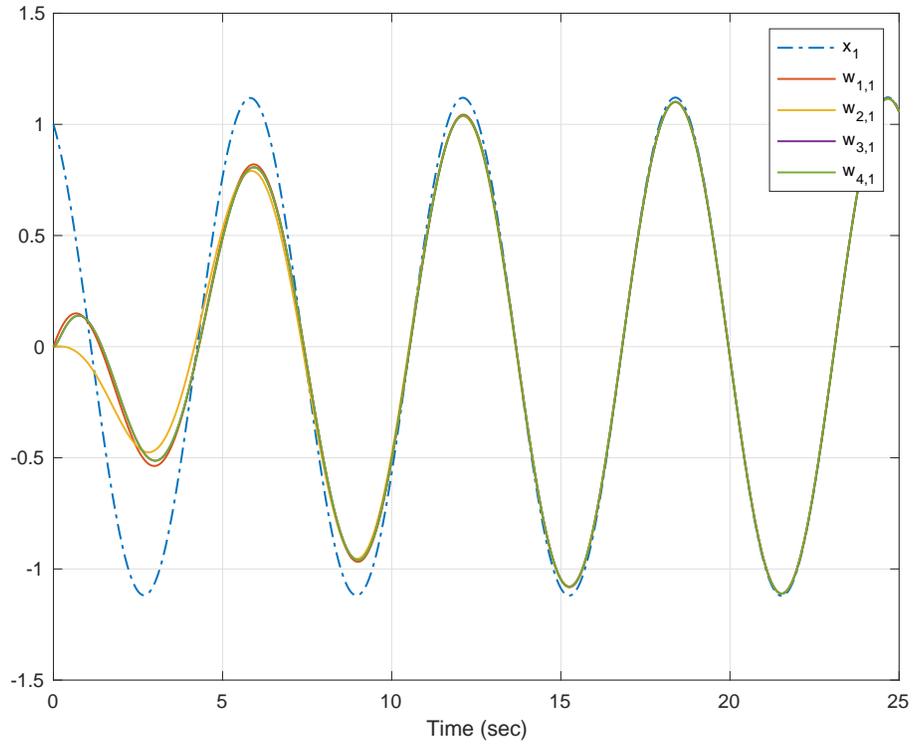}
	\caption{Plots of trajectories of $x_1$  (dashed lines) and the corresponding filter state component (solid lines).}
	\label{state1}
\end{figure}

\begin{figure}[t!]
	\centering
	\includegraphics[width=0.8\columnwidth]{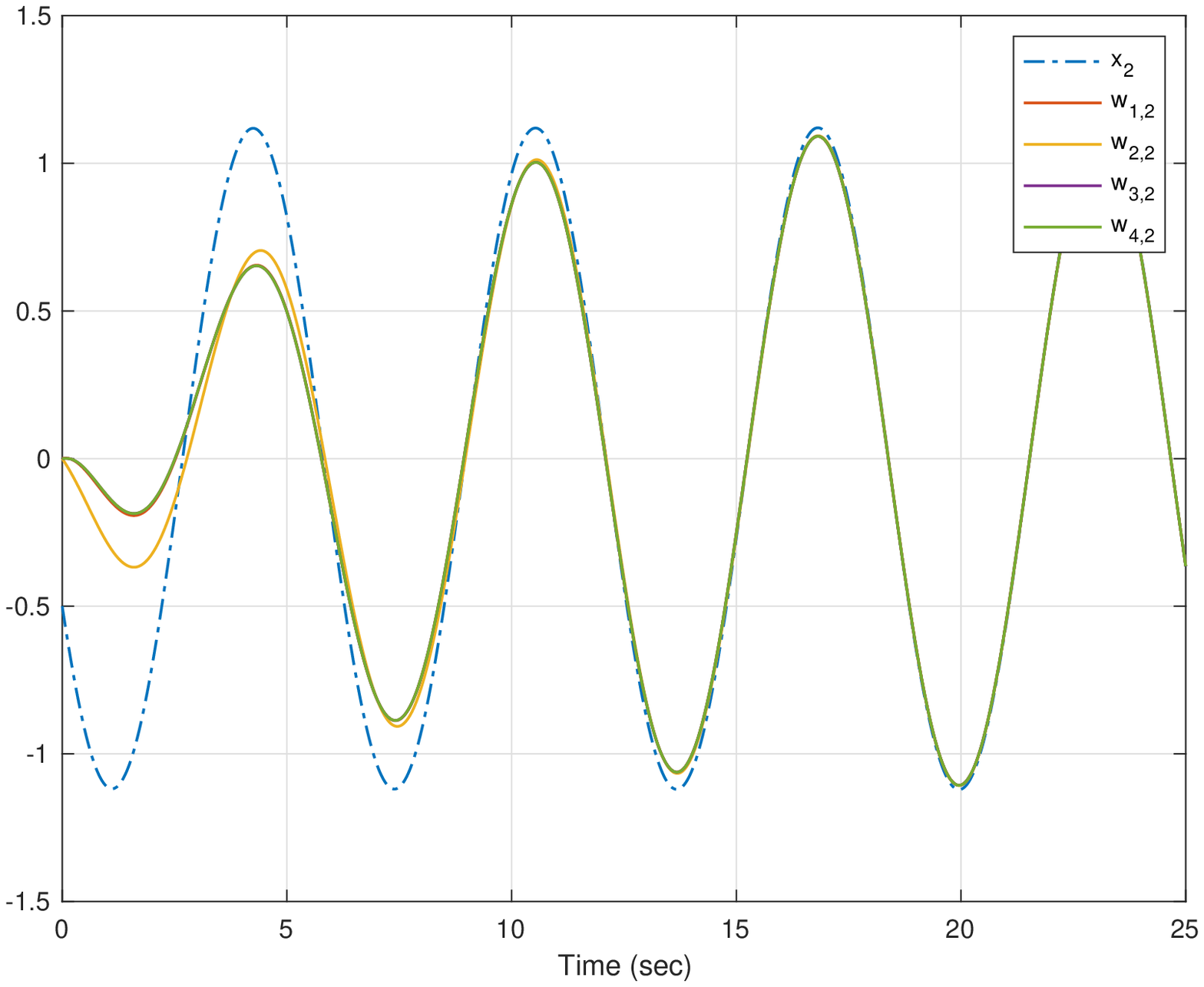}
	\caption{Plots of trajectories of $x_2$  (dashed lines) and the corresponding filter state component (solid lines).}
	\label{state2}
\end{figure}

\begin{figure}[t]
	\centering
	\includegraphics[width=0.8\columnwidth]{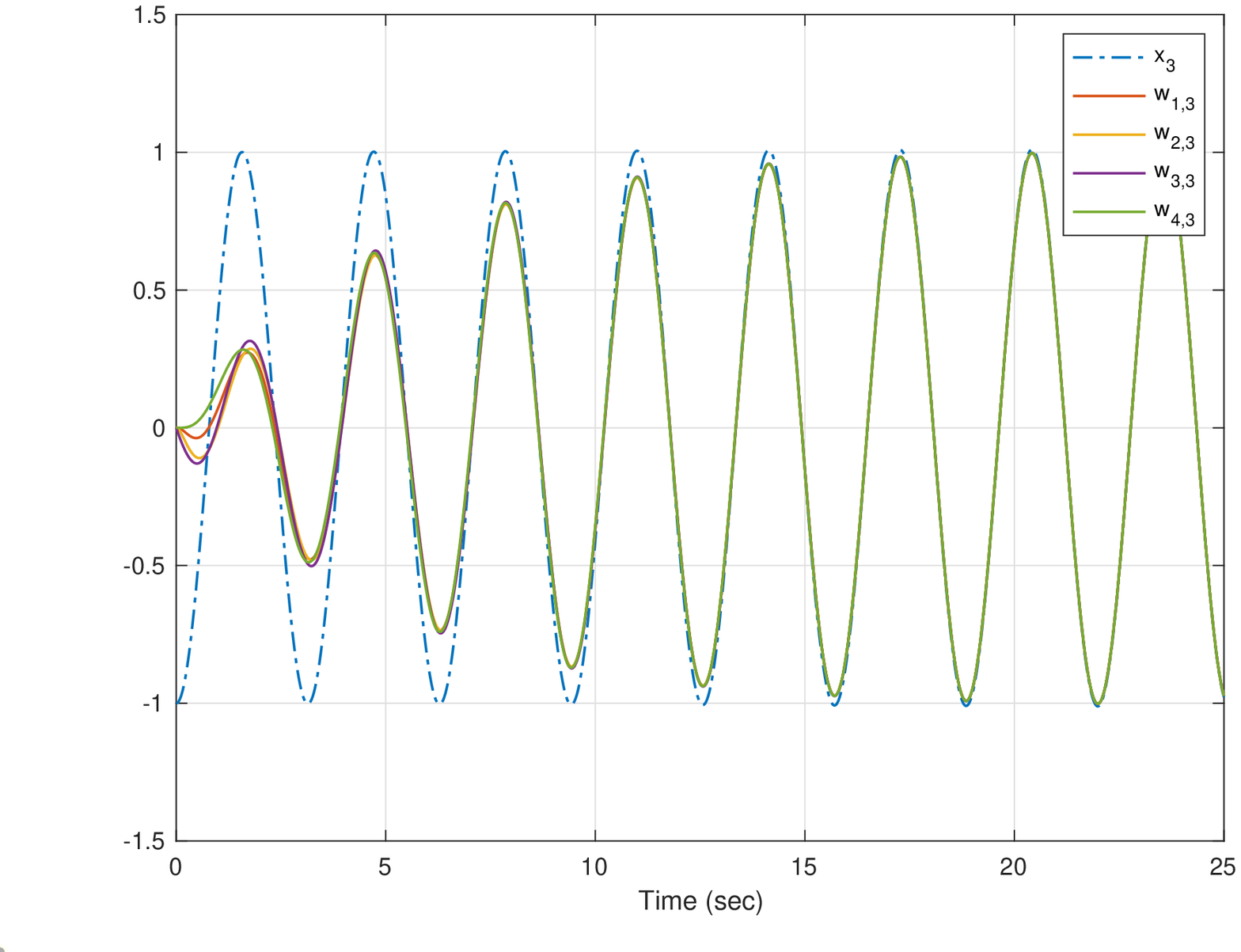}
	\caption{Plots of trajectories of $x_3$  (dashed lines) and the corresponding filter state component (solid lines).}
	\label{state3}
\end{figure}

\begin{figure}[t]
	\centering
	\includegraphics[width=0.8\columnwidth]{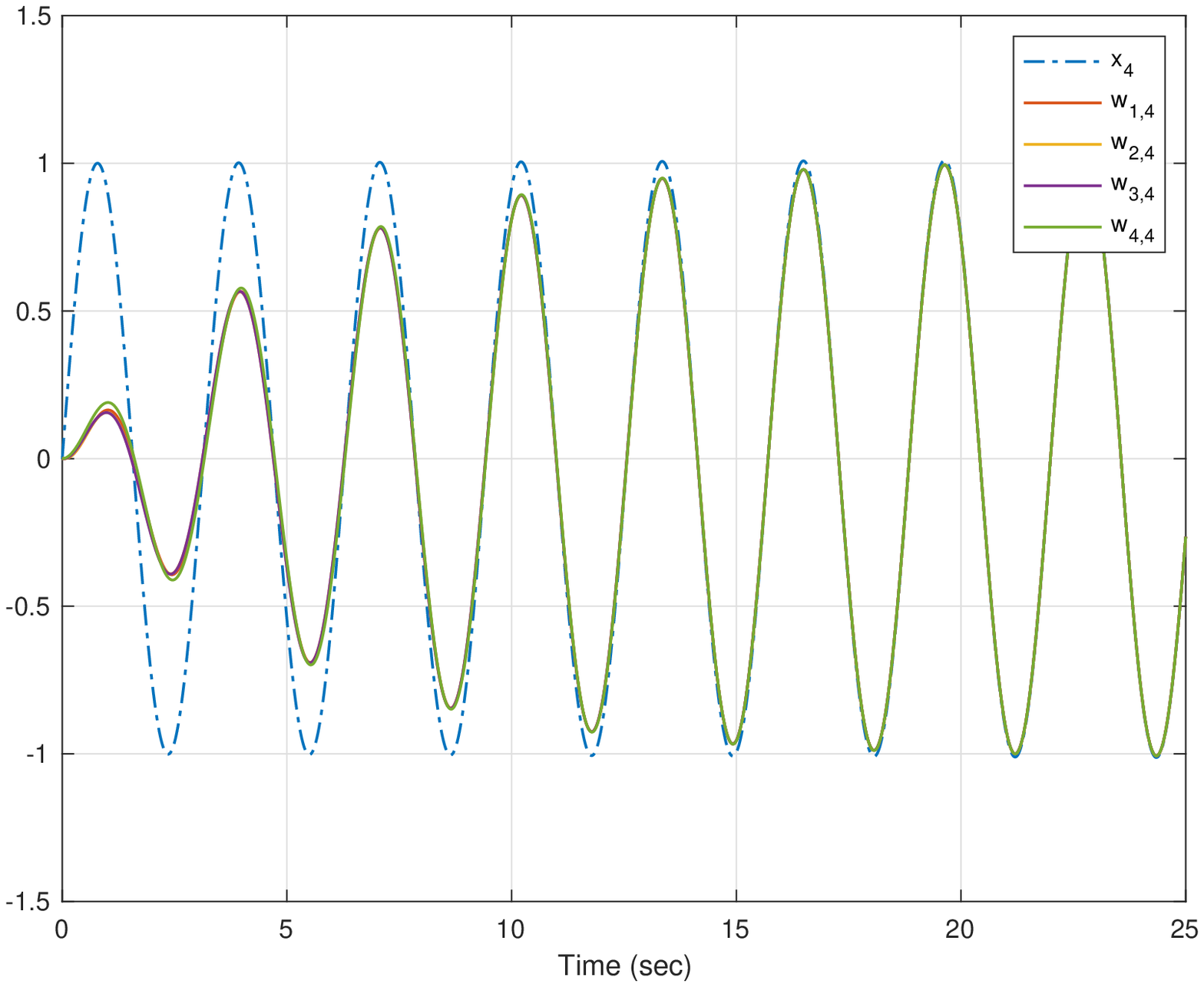}
	\caption{Plots of trajectories of $x_4$  (dashed lines) and the corresponding filter state component (solid lines).}
	\label{state4}
\end{figure}
As an example, we take the initial state of the system \eqref{sys_sim} to be $x_0 = \begin{pmatrix}
1 & -0.5 & -1 & 0
\end{pmatrix}^\top$ and  the initial state  of the distributed filter to be zero.
In Figures \ref{state1}, \ref{state2}, \ref{state3} and \ref{state4}, we have plotted the state trajectories of the system and that of the distributed filter in  absence of  external disturbances.
It can be seen that the states of the local filters asymptotically track  the state of the system  \eqref{sys_sim}.
Moreover, we compute 
\begin{equation*}
\begin{aligned}
\sum_{i=1}^{N}\textnormal{tr} \left[ (E_{i1}^\top - D_i^\top  G_{i1}^\top) P_{i1} (E_{i1} - G_{i1}D_i)  + E_{i2}^\top  E_{i2} \right] =   1.3717.
\end{aligned}
\end{equation*}
Thus, for all $\gamma> 1.3717$, the distributed filter \eqref{dis_observer} with gain matrices $F_i$ and $G_i$ is $\mathcal{H}_2$ $\gamma$-suboptimal.


\section{Conclusion}\label{sec_conclusion}
In this paper, we have studied the $\mathcal{H}_2$ and $\mathcal{H}_\infty$ suboptimal distributed filtering problem for linear systems.
We have established conditions for the existence of suitable filter gains. These are expressed in terms of solvability of LMI's.
Based on these conditions, we have provided  conceptual algorithms for obtaining the $\mathcal{H}_2$ and $\mathcal{H}_\infty$ suboptimal distributed filters, respectively.
The computation of these distributed filters requires centralized computation, i.e. global information is needed. 
As a possibility for future research, we mention the extension of the results in this paper to the case that  the filter gains need  to be computed in a decentralized fashion, see for example \cite{8943159}.


%


%

\balance
%
%
%

\ifCLASSOPTIONcaptionsoff
  \newpage
\fi



%
%
%

\bibliographystyle{IEEEtran}
\bibliography{h2_observer}

%

%
%
%




\end{document}